\documentclass[12pt,twoside]{amsart}\usepackage{amssymb}
\usepackage{enumerate}

\newtheorem{thmspec}{\relax}

\newtheorem{theorem}{Theorem}[section]
\newtheorem{thm}[theorem]{Theorem}

\newtheorem{prop}[theorem]{Proposition}
\newtheorem{defi}[theorem]{Definition}

\theoremstyle{definition}

\theoremstyle{remark}

\numberwithin{equation}{section}
\def \Bbb{\mathbb}

\def\onto{{\kern3pt\to\kern-8pt\to\kern3pt}}

\def\<{\langle}
\def\>{\rangle}
\def\|{{\ |\ }}

\def\onto{\twoheadrightarrow}

\def\-{\underline}

\def\N{\Bbb N}

\def\R{\Bbb R}
\def\C{\Bbb C}

\def\d{{\rm d}}

\def\P{\mathbb{P}}

\def\K{\mathcal{K}}

\def\<{\langle}
\def\>{\rangle}

\catcode`\@=11
\def\serieslogo@{\relax}
\def\@setcopyright{\relax}
\catcode`\@=12

\title[Thermodynamics of Rational Maps]
{On Thermodynamics of Rational Maps on the Riemann sphere}

\date{February 20, 2006}

\begin{document}

\author{Tien-Cuong  Dinh}
\address{Tien-Cuong  Dinh \\
  Analyse Complexe, Plateau 7D, Institut de Math{\'e}matique de Jussieu,
  175 rue du Chevaleret, 75013 Paris, France}
\email{dinh@math.jussieu.fr}

\author{Vi{\^e}t-Anh  Nguy{\^e}n}
\address{Vi{\^e}t-Anh  Nguy{\^e}n\\
 Max-Planck Institut f{\"u}r    Mathematik\\
Vivatsgasse 7, D--53111\\
Bonn, Germany}
\email{vietanh@mpim-bonn.mpg.de}

\author{Nessim Sibony}
\address{Nessim Sibony \\
  Math{\'e}matique-B{\^a}t 425\\
 UMR 8628, Universit{\'e} Paris-Sud
  F--91405 Orsay, France}
\email{nessim.sibony@math.u-psud.fr}

\subjclass[2000]{Primary 37A25, 37F, Secondary  32U40, 32H50}

\keywords{(Local) Central Limit Theorem, spectral gap}

\maketitle


\begin{abstract}
We show that there is a spectral gap for the transfer operator
associated to a rational map $f$ on the Riemann sphere.
Using this and the method of pertubed operators  we  establish the   
(Local) Central Limit Theorem for the measure of maximal entropy of
$f$, with an estimate on the  speed of convergence.
 \end{abstract}


\section{Introduction}

Let $f:\P^1\rightarrow\P^1$ 
be a rational map of degree $d\geq 2$. The Riemann sphere
$\P^1$ is endowed  with the  Fubini-Study form  $\omega$, normalized by
$\int_{\P^1}\omega=1$. It is well-known that the sequence
$d^{-n}f^{n*}(\omega)$ converges weakly to
a probability measure $\mu$ which is {\it totally
  invariant} : $d^{-1}f^*(\mu)=f_*\mu=\mu$, mixing, and which maximizes the
entropy. This is the {\it equilibrium measure} of $f.$
See Lyubich \cite{ly} or \cite{si, ds2, ds1}  for the multidimensional analogue.

 The classical Birkhoff ergodic theorem
 implies that if $\phi$ is a $\mu$-integrable function then the Birkhoff sum 
$ \frac{1}{n}\sum_{i=0}^{n-1}\phi\circ f^i$
 converges $\mu$-almost everywhere ($\mu$-a.e. for short) to $\mu(\phi)$. 
Assume that   $\mu(\phi)=0.$ Then,
the Birkhoff sum converges $\mu$-a.e to $0.$

Recall that $\phi$ is a {\it coboundary}  if there exists
  a   function $\psi\in L^2(\mu)$ such that $\phi=\psi\circ f-\psi.$
  One  easily checks that  for such a function 
$ \frac{1}{\sqrt{n}}\sum_{i=0}^{n-1}\phi\circ f^i$ converges to $0$  $\mu$-almost surely.
  We say that  a real-valued function $\phi \in L^1(\mu)$ satisfies
  the {\it Central Limit Theorem} 
(CLT for short) if 
   $ \frac{1}{\sqrt{n}}\sum_{i=0}^{n-1}\phi\circ f^i$  
converges in law to a Gaussian random variable 
   of zero mean and variance $\sigma>0.$  That is, for every interval $I\subset \R,$
   \begin{equation*}
   \lim\limits_{n\to\infty} \mu \Big\lbrace
   \frac{1}{\sqrt{n}}\sum\limits_{i=0}^{n-1}\phi\circ f^i\in I  \Big\rbrace
   =\frac{1}{\sqrt{2\pi}\sigma}\int\limits_{I}e^{-\frac{t^2}{2\sigma^2}}
   \d  t.
   \end{equation*}
   Such a function is not a coboundary.
   
Let $\Lambda$ denote the  adjoint of $f^*$ acting on $L^2(\mu)$.
Sometimes, $\Lambda$ is  called the {\it Perron-Frobenius operator}. We have
   \begin{equation} \label{eq_transfert}
   \Lambda\phi(z)=\frac{1}{d}(f_{\ast}\phi)(z)=\frac{1}{d}\sum\limits_{f(w)=z} \phi(w),\qquad z\in\P^1.
 \end{equation} 
   Consider  a real-valued function $\phi\in L^{\infty}(\mu)$ such
   that  $\mu(\phi)=0$, which is not a coboundary. A classical theorem of
   Gordin-Liverani \cite{go, li} 
   implies that if
   \begin{equation} \label{eq_gordin}
   \sum\limits_{n\geq 0}\Vert \Lambda^n\phi\Vert_{L^1(\mu)}<\infty, 
   \end{equation} 
  then $\phi$ satisfies the CLT with
  \begin{equation*}
  \sigma^2:=-\int\limits_{\P^1}\phi^2\d\mu+2\sum\limits_{n=0}^{\infty}\int\limits_{\P^1} 
\phi(\phi\circ f^n)\d \mu =\lim\limits_{n\to\infty}\int\limits_{\P^1}
 \Big( \frac{1}{\sqrt{n}} \sum_{k=0}^{n-1} \phi\circ f^k  \Big)^2 \d\mu .
  \end{equation*} 

When $\phi$ is smooth, condition (\ref{eq_gordin}) is
a direct consequence of the exponential decay of correlations 
proved by Forn\ae ss and the third named author in
\cite{fs}.    
The case of H{\"o}lder continuous observables follows by interpolation between Banach
spaces \cite{do, ds1}, see also \cite{ds2, ds3} for the case of
meromorphic maps in higher dimension.
Note that with a completely different method (analysis of inverse
branches) Haydn \cite{ha} has proved the exponential decay of
correlations, hence the CLT, for Gibbs measures.
A weaker result was previously proved by Denker,
Przytycki and Urba\'nski \cite{de}.
   
In this work, we  study the speed of convergence in the  CLT 
and prove the Local CLT. 
The approach is  spectral theory for pertubation of the
Perron-Frobenius operator $\Lambda$ in appropriate Sobolev spaces.
We follow  the Nagaev method   \cite{na}, which was used by     
Rousseau-Eg{\`e}le \cite{re} and  Broise
\cite{br}.
The main point
   here is to prove some necessary estimates in order to apply the Nagaev
   method. 
It turns out that the analogous estimates in the 
   higher dimension case are more difficult. We also use a
   version of the classical result of Ionescu Tulcea and Marinescu.

   We close  this introduction with a brief outline of the paper to follow. 
In Section 2 we
collect some background on Sobolev spaces. 
In Section 3 we develop the method of spectral theory for pertubation of Perron-Frobenius operator.
This  will enable us to prove the CLT (Theorem A)
and   obtain  an estimate on the  speed of convergence  (Theorem B).
Section 4  develops auxiliary tools used in Section 5  to
  prove   the Local CLT  (Theorem C). 
   
\medskip
   
\noindent
{\bf Notation.} 
For a normed-space $E$, denote  by $\Vert\cdot\Vert_{E}$ its norm.
We write $L^p$ instead of $L^p(\P^1)$, where
the $L^p$-norm of a form is the sum of its coefficients for a fixed
atlas of $\P^1.$ If $\alpha_{10}$, $\alpha_{01}$ are forms of bidegrees
$(1,0)$ and $(0,1)$ respectively, and
$\alpha=\alpha_{10}+\alpha_{01}$, we use the equivalent $L^2$-norm
$\Vert\alpha\Vert_{L^2}=\int_{\P^1}
(i\alpha_{10}\wedge\overline\alpha_{10}+i\overline\alpha_{01}\wedge
\alpha_{01})$.
If $T:E\rightarrow E$ is a linear operator, we denote by $\rho(T)$ its
spectral radius. We also denote by $\Lambda$ the operator ${1\over d}
f_*$ even when it acts on forms. 
For a Radon measure $\nu$,  $\Vert\nu\Vert$ is its total variation.
The constants $C$, $C'$, $C''$ and $C_\epsilon$  are
not necessarily the same  at each occurence; they  depend  possibly on some parameters.


\section{Sobolev spaces}

The Sobolev space  $W^{1,p}$, $p\geq 1$, is the space of  complex-valued functions $\phi$ in $L^p$
such that  $\d\phi,$  which is defined in the sense of currents, has
$L^p$ coefficients, see \cite{af, he}. It is endowed with the
canonical norm $\Vert\phi\Vert_{W^{1,p}}=\Vert\phi\Vert_{L^p} + \Vert
\d\phi\Vert_{L^p}$. 
Lipschitz functions
belong to $W^{1,p}$.
Let $\mathcal{L}$ denote the Lebesgue measure of mass 1 on $\P^1$ associated to
the Fubini-Study form $\omega$. Define for each  $L^1$ function its {\it mean  value}
$m(\phi):= \int_{\P^1}\phi \d\mathcal{L}$.  We have the following Poincar{\'e}--Sobolev inequality:
for every real number $p\geq 1,$ there exists a constant $C=C(p)$ such that
\begin{equation}\label{PS_inequality}
\Big( \int\limits_{\P^1}\vert  \phi-m(\phi)\vert^p  \d\mathcal{L}   \Big)^{\frac{1}{p}}
\leq  C\Vert \d\phi\Vert_{L^2}, \qquad   \phi\in W^{1,2}.
\end{equation}
In particular,  
\begin{equation}\label{PS2_inequality}
\Vert \phi\Vert_{L^p}\leq |m(\phi)|+ C\Vert \d \phi\Vert_{L^2} 
\leq\Vert \phi\Vert_{L^1}+C\Vert \d \phi\Vert_{L^2}  \leq C'
\Vert\phi\Vert_{W^{1,2}}. 
\end{equation}
Hence, $W^{1,2}\subset L^p.$ Moreover,  by Sobolev embedding theorem, the latter inclusion is compact.
 For every $\epsilon>0,$  there is a constant  $C=C(\epsilon)$ such that
 \begin{equation}\label{PS3_inequality}
 \Vert \phi\Vert_{L^{\infty}}\leq C \Vert
 \phi\Vert_{W^{1,2+\epsilon}}, 
\qquad\phi\in  W^{1,2+\epsilon},
 \end{equation} 
 and 
 \begin{equation}\label{eq_Banach_algebra_W}
 \Vert \phi\psi\Vert_{W^{1,2+\epsilon}}\leq C\Vert 
\phi \Vert_{W^{1,2+\epsilon}}\Vert  \psi\Vert_{W^{1,2+\epsilon}},
 \qquad\phi,\psi\in W^{1,2+\epsilon},
 \end{equation}  
 (see, for example, Theorem 4.39 in \cite{af}).
 
We recall from Lemmas 5.3 and  5.4 in \cite{ds3} that
\begin{equation}\label{DS_inequality_W_1}
 \vert  \mu(\phi) \vert\leq  C\Vert \phi\Vert_{W^{1,2} }
\end{equation}  
and
\begin{equation}\label{DS_inequality_W_2}
\big\Vert\Lambda^n\phi -  \mu(\phi)\big\Vert_{W^{1,2} }\leq
\frac{C}{d^{n/2}}
\Vert \phi\Vert_{W^{1,2} }
\end{equation} 
for  every  continuous $\phi\in W^{1,2} $  and $n\geq 0.$  Since
smooth functions are dense in $W^{1,2}$ 
we can extend $\mu$ to a continuous linear form of $W^{1,2}$ in a
canonical way.
Therefore,  (\ref{DS_inequality_W_1})--(\ref{DS_inequality_W_2})
are still valid for any $\phi\in W^{1,2}$.
Moreover, $\mu$ is strongly continuous (WPC in the sense of
\cite{ds3}), i.e. if $\phi_n\rightarrow \phi$ weakly and
$\Vert\phi_n\Vert_{W^{1,2}}$ is bounded then
$\mu(\phi_n)\rightarrow\mu(\phi)$. It follows that one can associate to each
function in $W^{1,2}$ a function in $L^1(\mu)$ in a canonical way.

\begin{prop}\label{equivalent_norms}
The following norm 
\begin{equation*}
\Vert \phi\Vert :=
 \vert  \mu(\phi)\vert +  \Vert \d \phi\Vert_{L^{2}}
\end{equation*}
is equivalent to the standard norm of $W^{1,2}.$
\end{prop}
\begin{proof}  Inequality (\ref{DS_inequality_W_1}) implies that 
$\Vert\phi\Vert\leq C\Vert\phi\Vert_{W^{1,2}}$. From (\ref{PS2_inequality}), 
we only need to check that
there is a constant $C$ such that
\begin{equation}\label{eq_equivalent_norms}
\vert m( \phi)\vert \leq  
 \vert  \mu (\phi)\vert + C \Vert \d \phi\Vert_{L^{2}}, \qquad  \phi\in W^{1,2}.
\end{equation}
To prove this, observe that  $m(\phi-m(\phi))=0.$ Therefore, by (\ref{DS_inequality_W_1}),
\begin{equation*} 
\Big\vert \int\limits_{\P^1}   (\phi-m(\phi))\d\mu    \Big\vert
\leq  C\Vert \d\phi\Vert_{L^2},
\end{equation*}
which implies (\ref{eq_equivalent_norms}).
\end{proof}


\section{Central Limit  Theorem  and speed of convergence}

For every $\theta\in\C$ and  every real-valued function $h,$ consider the 
 pertubation of the Perron-Frobenius operator  $\Lambda$
\begin{equation}\label{pertubated_oper}
\Lambda_h(\theta)(\phi):= \Lambda(e^{\theta h} \phi).
\end{equation}

We collect here some properties of the family of   pertubations of $\Lambda.$

\begin{prop}\label{prop_pertubated_oper}
Assume that  $h$  is in $W^{1,2+\epsilon}$ 
for some $\epsilon>0.$  Then
\begin{enumerate}
\item[(1)]   $\Lambda_h(0)=\Lambda.$
\item[(2)] For every $\theta\in\C,$  $\Lambda_h(\theta)$ is  a bounded linear operator on $W^{1,2}$.
\item[(3)] The maps $\theta\mapsto  \Lambda_h(\theta)$ is  analytic.
\item[(4)] For every $n\in\N,$   $\theta\in\C$  and every $\phi$ in $W^{1,2}, $ 
  we have
   \begin{equation*}
   \Lambda_h^n(\theta)(\phi)= \Lambda^n(e^{\theta S_n h} \phi),
   \end{equation*}
where $S_nh:= \sum\limits_{k=0}^{n-1} h\circ f^k$ for $n\geq 1$ and $S_0h:=0$.
\end{enumerate}
\end{prop}
\begin{proof} For the proof of Parts (2) and (3) we will need the following inequality
\begin{equation}\label{eq1_prop_pertubated_oper}
\Vert \phi\psi\Vert_{W^{1,2 } }\leq C_{\epsilon}\Vert \phi \Vert_{W^{1,2 } }\Vert  \psi\Vert_{W^{1,2+\epsilon} },
 \qquad\phi\in W^{1,2} ,\psi\in W^{1,2+\epsilon} .  
\end{equation}
This  inequality is a consequence of the following  estimates where we
use (\ref{PS2_inequality}), (\ref{PS3_inequality}) and the classical
H{\"o}lder inequality
\begin{eqnarray*}
\Vert \phi\psi\Vert_{L^{2 } }&\leq &  \Vert \phi \Vert_{L^{2 } }\Vert  \psi\Vert_{L^{\infty} }
\leq  C_{\epsilon}\Vert \phi \Vert_{L^{2 } }\Vert
\psi\Vert_{W^{1,2+\epsilon} }
\end{eqnarray*}
and
\begin{eqnarray*}
\Vert \d(\phi\psi)\Vert_{L^{2 }  }&\leq &  \Vert\psi\d \phi
\Vert_{L^{2 } }+\Vert  \phi\d\psi\Vert_{L^{2} }\\
&\leq&  C_{\epsilon}\Vert  \psi\Vert_{L^{\infty} }\Vert\d \phi \Vert_{L^{2 } }+
\Vert\d  \psi\Vert_{L^{2+\epsilon} }\Vert \phi\Vert_{L^{\frac{2(2+\epsilon)}{\epsilon}} },\\
&\leq & C_{\epsilon}\Vert  \psi\Vert_{ W^{1,2+\epsilon} }\Vert\d \phi \Vert_{L^{2 } }+C_{\epsilon}
\Vert\d  \psi\Vert_{L^{2+\epsilon} }\Vert \phi\Vert_{ W^{ 1,2} }.
\end{eqnarray*}

 Using  (\ref{eq1_prop_pertubated_oper})and  (\ref{eq_Banach_algebra_W}),
we get 
$$\Vert  \Lambda(h^n \phi)\Vert_{W^{1,2} }
\leq  \Vert \Lambda\Vert_{W^{1,2} }  \Vert h^n\phi\Vert_{W^{1,2} }
\leq  C\Vert \Lambda\Vert_{W^{1,2} } 
\Vert  \phi\Vert_{W^{1,2} } \left( C \Vert 
h\Vert_{W^{1,2+\epsilon} }\right)^n. $$
Consequently, the series $\sum_{n\geq 0}\frac{\theta^n}{n!}\Lambda(h^n \phi)$ 
converges normally in $W^{1,2}$. Its limit is
equal to  $\Lambda_h(\theta)(\phi),$  and
\begin{equation*}
\Vert   \Lambda_h(\theta)(\phi)\Vert_{W^{1,2} } \leq
C e^{C\vert \theta\vert \!\!\ \Vert h\Vert_{W^{1,2+\epsilon} }}
\Vert\Lambda\Vert_{W^{1,2}} \Vert\phi\Vert_{W^{1,2} }.
\end{equation*}
This gives (2) and (3).

Part (4)   is proved  by induction. Using (\ref{eq_transfert}), we get
$\Lambda^n\big((h\circ f^n)\phi\big)= h \Lambda^n\phi$ for $n\geq 0$.
Hence
\begin{eqnarray*}
\Lambda^n(e^{\theta S_n h}\phi) & = & \Lambda\big(\Lambda^{n-1}(e^{\theta
  h\circ f^{n-1}} e^{\theta S_{n-1} h}\phi)\big)\\
& = &   \Lambda\big( e^{\theta h}\Lambda^{n-1} (e^{\theta S_{n-1} h}\phi)\big) \\
& = & \Lambda_h(\theta)\big(    \Lambda^{n-1}(  e^{\theta S_{n-1}
  h}\phi)   \big).
\end{eqnarray*}
This completes the proof.
\end{proof}

The spectrum of $\Lambda_h(\theta)$ is described  by the following pertubation theorem. 

\begin{thm}\label{thm_pertubation}
Let $h$   be a
real-valued function in $W^{1,2+\epsilon}. $  
  Then there exists a real number $a>0$  such that  for 
$\vert\theta\vert <a $  and $n\in \N,$ we have the following decomposition
\begin{equation*}
\Lambda_h^n(\theta)=\lambda^n(\theta)\Phi(\theta) +\Psi^n(\theta) ,
\end{equation*}
where
\begin{itemize}
\item  
$\lambda(\theta)$ is the eigenvalue with maximal modulus  of
$\Lambda_h(\theta),$ $\lambda(0)=1,$    
and  $\vert  \lambda(\theta)\vert>
\frac{2+\rho(\Psi(0))}{3}=:\rho_1 <1;$
\item  $\Phi(\theta)$  is the projection of $W^{1,2}$  onto the eigenspace $W^{1,2}(\lambda(\theta))$
which is of dimension $1$, and $\Phi(0)(\phi)=\mu(\phi);$
\item $\Psi(\theta)$ is a bounded linear operator on $W^{1,2}$ with $\rho(\Psi(\theta))< \rho_1,$  and
\begin{equation*}
\Psi(\theta)\circ\Phi(\theta)=\Phi(\theta)\circ\Psi(\theta)=0.
\end{equation*}
\item   $\theta\mapsto \Phi(\theta),  $ $\theta\mapsto \Psi(\theta),  $
$\theta\mapsto \lambda(\theta)$  are  analytic maps on the disc $\{\vert\theta\vert <a\}.$
\end{itemize}
 \end{thm}
 \begin{proof} 
Observe that $\Lambda(1)=1$. Then, constant functions are
eigenfunctions associated to the eigenvalue 1 of the operator
$\Lambda_h(0)=\Lambda$. If $\lambda(0)$, $\Phi(0)$ are
as above and if $\Psi(0):=\Lambda-\Phi(0)$, then
 \begin{equation*}
 \Lambda_h(0)=\Lambda=\lambda(0) \Phi(0)+\Psi(0),
 \end{equation*}
It follows from
   (\ref{DS_inequality_W_2})  
 that  $\rho(\Psi(0))\leq \frac{1}{\sqrt{d}}<1.$
Hence, $\lambda(0)=1$ is the unique eigenvalue of maximal modulus of
$\Lambda$ and $W^{1,2}(\lambda(0))$ is of dimension 1.
Using Proposition   \ref{prop_pertubated_oper}, it is enough to apply
   the Rellich pertubation method
described in \cite[Chapter VII]{dus}  (see also \cite[Proposition 5.2]{br}, \cite{lp}).
\end{proof}

 Now we are in  position to   prove the first main result.
  \renewcommand{\thethmspec}{Theorem A}
  \begin{thmspec} 
 Let $f:\P^1\rightarrow\P^1$ 
be a rational map of degree $d\geq 2$. Let $g$  be a real-valued function in $W^{1,2 } $
 and    $h$  a real-valued function in $W^{1,2+\epsilon}$,
 $\epsilon>0$, 
 such that  $\mu(g)=1$, $\mu(h)=0$ and   $g(z)\geq 0$ for $\mu$-a.e.
 $z\in \P^1$. Assume  
$h$ is not a coboundary and let $\sigma$ be the positive number defined by
 \begin{equation} \label{eq_sigma}
 \sigma^2:=-\int\limits_{\P^1}h^2\d \mu+ 2\sum\limits_{n\geq 0}
 \int\limits_{\P^1} h(h\circ f^n) \d \mu=\lim\limits_{n\to\infty}\int\limits_{\P^1}
 \Big( \frac{1}{\sqrt{n}} S_nh  \Big)^2 \d\mu.
 \end{equation}  
 Define the   probability  measure
 \begin{equation*}
 \nu:=g\mu.
 \end{equation*} 
 Then for every real number $v,$  we have
 \begin{equation*}
 \lim\limits_{n\to\infty}\nu\Big\lbrace  \frac{S_n h }{\sigma\sqrt{n}}\leq v \Big\rbrace
 =\frac{1}{\sqrt{2\pi}}\int\limits_{-\infty}^v  e^{-\frac{t^2}{2}} \d t.
 \end{equation*}
 \end{thmspec}

When $g=1$, we obtain the CLT for $\mu$. Since the theorem
holds for every $\nu=g\mu$, the set $\lbrace  S_n h  \leq \sigma\sqrt{n}v \rbrace$
is asymptotically equidistributed
with respect to $\mu$ as $n\rightarrow\infty$, because this set is
independent of $g$.

 \begin{proof}
 Using Proposition  \ref{prop_pertubated_oper} and Theorem \ref{thm_pertubation}
 the proof follows  the same lines as  in \cite[Th{\'e}or{\`e}me 6.8]{br} or
 \cite[Th{\'e}or{\`e}me 2]{re}.  For the sake of clarity,  we  give here
the main aguments. The proof is divided into four steps.
 
\medskip

\noindent
 {\bf Step 1:}  Proof of the equality $\lambda^{'}(0)=0.$
Since $\mu$ is totally invariant, we have $\int\Lambda^n(\phi)\d \mu
=\int \phi\d \mu$ for every $\phi\in L^1(\mu)$. This, 
Proposition \ref{prop_pertubated_oper} and 
Theorem \ref{thm_pertubation} imply that for $\left\vert \frac{t}{n}\right\vert <a,$ 
\begin{eqnarray*}
\int_{\P^1} e^{itS_nh\over n}\d \mu & = & \int_{\P^1} \Lambda^n\big( e^{itS_nh\over
  n} \big)\d \mu = \int\limits_{\P^1}\Lambda_h^n\Big( \frac{it}{
  n}\Big)(1)\d\mu \\
& = & \lambda^n\Big( \frac{it}{n}\Big) 
\int\limits_{\P^1}\Phi\Big( \frac{it}{ n}\Big)(1)\d\mu  +
\int\limits_{\P^1}\Psi^n \Big( \frac{it}{ n}\Big)(1)\d\mu .
\end{eqnarray*}
We then use  the Taylor expansion of order $2$ at $0$  of
$\lambda\big( \frac{it}{ n}\big)  $  
and $\Phi\big( \frac{it}{ n}\big).$
Moreover, we use the estimate  
$ \Vert \Psi^n \big( \frac{it}{
      n}\big)(\phi)\Vert_{W^{1,2}}\leq C\rho_1^n.$
   Therefore, one obtains that
  \begin{equation*}
  \lim\limits_{n\to\infty} \int\limits_{\P^1}e^{it \frac{S_nh}{n}} \d\mu=e^{it\lambda^{'}(0)}.
  \end{equation*}
  By Birkhoff's theorem, $\lim_{n\to\infty} \frac{S_nh}{n}=\mu(h)=0,$  $\mu$-a.e. 
 Hence,  $\lambda^{'}(0)=0.$

\medskip

\noindent
 {\bf Step 2:}  Proof of the equality $\lambda^{''}(0)=\sigma^2.$
Equality (\ref{eq_sigma}) is easy to check using the invariance of
$\mu$.
 On the other hand,
 \begin{equation}\label{eq2_Step2_TheoremA}
 \int\limits_{\P^1}\Big(  \frac{S_n h}{\sqrt{n}}\Big)^2\d\mu
=-\frac{\partial^2}{\partial t^2}\Big\lbrack\int\limits_{\P^1}
 e^{\frac{it}{\sqrt{n}} S_nh }  \d\mu\Big\rbrack_{t=0}.
 \end{equation}
 By Theorem \ref{thm_pertubation}, we may rewrite the expression in
 brackets 
in (\ref{eq2_Step2_TheoremA}) as follows
 \begin{eqnarray*}
 \int\limits_{\P^1}
 e^{\frac{it}{\sqrt{n}} S_nh }  \d\mu & = & 
\int\limits_{\P^1}\Lambda_h^n\Big( \frac{it}{
 \sqrt{n}}\Big)(1)\d\mu \\
& = & \lambda^n\Big( \frac{it}{\sqrt{n}}\Big) 
\int\limits_{\P^1}\Phi\Big( \frac{it}{ \sqrt{n}}\Big)(1)\d\mu  
+\int\limits_{\P^1}\Psi^n \Big( \frac{it}{ \sqrt{n}}\Big)(1)\d\mu .
\end{eqnarray*}
The derivative of the second term 
 can be bounded by $C\rho_1^n$ using Cauchy's formula.
We use the Taylor expansion  of order $2$ at $0$  of  
$\lambda\big( \frac{it}{\sqrt{ n}}\big)$  and $\Phi\big( \frac{it}{\sqrt{ n}}\big)$
 and insert them into the above expansion.  After   taking the second derivative     we see that
 the right hand   side of (\ref{eq2_Step2_TheoremA}) is equal to
 $\lambda^{''}(0).$ We use here the equality $\lambda'(0)=0$.
  
\medskip

\noindent
 {\bf Step 3:}  Proof of the equality $ \lim_{n\to\infty} 
\int \Lambda_h^n\big( \frac{it}{ \sqrt{n}}\big)(g)\d\mu 
=e^{-\frac{\sigma^2 t^2}{2}}.$
By Theorem \ref{thm_pertubation}, we  have 
 \begin{equation}\label{eq1_Step3_TheoremA}
  \int\limits_{\P^1}\Lambda_h^n\Big( \frac{it}{
    \sqrt{n}}\Big)(g)\d\mu=
\lambda^n\Big( \frac{it}{\sqrt{n}}\Big) 
\int\limits_{\P^1}\Phi\Big( \frac{it}{ \sqrt{n}}\Big)(g)\d\mu  
+\int\limits_{\P^1}\Psi^n \Big( \frac{it}{ \sqrt{n}}\Big)(g)\d\mu .
\end{equation}
 As in Step 1, the last term in the above expression tends to $0$ as $n\to\infty.$ On the other hand,
using the previous steps and the hypothesis $\mu(g)=1$,  a straightforward computation gives that
$$  \lambda^n\Big( \frac{it}{\sqrt{n}}\Big)  = 
e^{-\frac{\sigma^2 t^2}{2}\big\lbrack 1 + O( \frac{t}{ \sqrt{n}})
 \big\rbrack}
\quad \mbox{and}
\quad
  \int\limits_{\P^1}\Phi\Big( \frac{it}{ \sqrt{n}}\Big)(g)\d\mu
=1+ O\Big( \frac{t}{ \sqrt{n}}\Big),$$
Replacing them in the right hand side of 
(\ref{eq1_Step3_TheoremA}), we obtain the desired equality.

\medskip

\noindent
{\bf Step 4:}  Conclusion. Using Proposition \ref{prop_pertubated_oper}(4), 
 we may rephrase the result of Step 3 as
   $ \lim_{n\to\infty} \int  e^{it \frac{S_nh}{
  \sigma\sqrt{n}}} d\nu 
=e^{-\frac{ t^2}{2}}.$
   Since the right hand side of this equality is the Fourier transform
  of the normal law  
$\mathcal{N}(0,1),$
   the Paul L{\'e}vy's method of characteristic functions implies that the
  sequence of 
random variables $\big( \frac{S_nh}{ \sigma\sqrt{n}}\big)_{n=1}^{\infty}$
   converges in law to the normal law. This implies the result.
  \end{proof}

 We come to the second main result.
  \renewcommand{\thethmspec}{Theorem B}
  \begin{thmspec} 
  Under the hypothesis of Theorem A, there is a  constant $C$  such that
 \begin{equation*}
 \sup\limits_{v\in\R}\Big\vert\nu\Big\lbrace  \frac{S_n h}{\sigma\sqrt{n}}\leq v \Big\rbrace
 -\frac{1}{\sqrt{2\pi}}\int\limits_{-\infty}^v  e^{-\frac{t^2}{2}}
 \d t\Big\vert \leq \frac{C}{\sqrt{n}}.
 \end{equation*}
 \end{thmspec}
 \begin{proof}
 First we recall the   Berry-Essen inequality (see \cite{fe}) in our
 context.  
Namely, there is a positive  constant $K$
 such that for $V>0,$
 \begin{equation*}
 \sup\limits_{v\in\R}\Big\vert\nu\Big\lbrace  \frac{S_n h }{\sigma\sqrt{n}}\leq v \Big\rbrace
 -\frac{1}{\sqrt{2\pi}}\int\limits_{-\infty}^v  e^{\frac{-t^2}{2}} \d t\Big\vert \leq 
 \frac{K}{V}+\frac{1}{\pi} \int\limits_{-V}^V 
\frac{1}{|t|}\!\ \Big\vert\int\limits_{\P^1} e^{it\frac{S_nh}{\sigma\sqrt{n}}} \d\nu  
 -e^{-\frac{t^2}{2}}  \Big\vert \d t .
 \end{equation*}

Choose $V=\sigma\sqrt{n}a'$ with $a'<a$ small enough, 
where  $a$ is the constant given by  Theorem \ref{thm_pertubation}.
 Then for every $t\in [-V,V],$  we have  $\big\vert\frac{t}{\sigma\sqrt{n}}\big\vert\leq
 a'<a.$ Consequently, applying Proposition
 \ref{prop_pertubated_oper}(4) and  Theorem  \ref{thm_pertubation} yields that
 \begin{eqnarray*}
\int\limits_{\P^1} e^{it\frac{S_nh}{\sigma\sqrt{n}}}\d\nu&=&\int\limits_{\P^1} e^{it\frac{S_nh}{\sigma\sqrt{n}}}g\d\mu= \int\limits_{\P^1}\Lambda_h^n\Big(\frac{it}{\sigma\sqrt{n}} 
\Big )(g)\d\mu\\
&=&\lambda^n\Big(   \frac{it}{\sigma\sqrt{n}}   \Big )\int\limits_{\P^1}\Phi\Big(\frac{it}{\sigma\sqrt{n}} \Big )(g)\d\mu  +\int\limits_{\P^1}\Psi^n\Big(\frac{it}{\sigma\sqrt{n}}
 \Big )(g)\d\mu.
 \end{eqnarray*}
Theorem  \ref{thm_pertubation} also implies that
$\int\Psi(0)(\phi)\d\mu =0$ for every $\phi$. Moreover,
the $W^{1,2}$-norm of the operator $\Psi(\theta)-\Psi(0)$ satisfies
$\Vert \Psi(\theta)-\Psi(0)\Vert_{W^{1,2}} \leq
C\vert\theta\vert$. Hence
we can apply this to $\phi:=\Psi^{n-1}\big(\frac{it}{\sigma\sqrt{n}}
 \big )(g)$ and obtain
   \begin{equation*}
 \Big\vert\int\limits_{\P^1}\Psi^n\big(\frac{it}{\sigma\sqrt{n}}
 \big )(g)\d\mu \Big\vert
=  \Big\vert\int\limits_{\P^1} \Big[\Psi\big(\frac{it}{\sigma\sqrt{n}}
 \big )-\Psi(0)\Big] \circ\Psi^{n-1}\big(\frac{it}{\sigma\sqrt{n}}
 \big )(g)\d\mu \Big\vert
\leq \frac{C\rho_1^{n-1}|t|}{\sqrt{n}}.
 \end{equation*}
Now, we use the computation in the third step  of the proof of Theorem
A. From the inequalities $1-e^{-sx}\leq e^{sx} -1\leq se^x$ for
$0\leq s\leq 1$ and $x\geq 0$, we deduce that
$$\Big \vert \lambda^n  \Big(   \frac{it}{\sigma\sqrt{n}}   \Big ) 
  -e^{-\frac{t^2}{2}}   \Big\vert \leq e^{-{t^2\over 2}}\Big\vert 
    e^{O\big({t^3\over \sigma\sqrt{n}}\big)}-1\Big\vert \leq
e^{-{t^2\over 2}}\Big(
    e^{{t^2\over 4} {C\vert t\vert\over \sigma\sqrt{n}}}-1\Big) \leq
    {Ce^{-{t^2\over 4}}\vert t \vert \over\sigma\sqrt{n}}$$
since ${C\vert t\vert\over\sigma\sqrt{n}}\leq Ca'\leq 1$ when
$a'$ is small enough.
It follows that
$$ \Big \vert \lambda^n  \Big(   \frac{it}{\sigma\sqrt{n}}   \Big ) 
  \int\limits_{\P^1}\Phi\Big(\frac{it}{\sigma\sqrt{n}} \Big )(g)\d\mu
  -e^{-\frac{t^2}{2}}   \Big\vert \leq \frac{Ce^{-{t^2\over 4}}|t|}{\sqrt{n}}.$$
  Hence,
   \begin{equation*}
   \int\limits_{-V}^V \frac{1}{|t|} \ \Big\vert\int\limits_{\P^1} 
e^{it\frac{S_nh}{\sigma\sqrt{n}}} \d\nu  
 -e^{-\frac{t^2}{2}}  \Big\vert \d t\leq \frac{C}{\sqrt{n}}
 \int\limits_{-V}^V 
 \big ( e^{-\frac{t^2}{4}}+ \rho_1^{n-1}\big) \d t
\leq \frac{C'(1+V\rho_1^{n-1})}{\sqrt{n}}.
 \end{equation*}
 Now we  substitute  $V=\sigma \sqrt{n} a'$ into the latter
 estimate. 
Consequently, it follows from the Berry-Essen inequality that
  \begin{equation*}
 \sup\limits_{v\in\R}\Big\vert\nu\Big\lbrace  \frac{S_n h}{\sigma\sqrt{n}}\leq v \Big\rbrace
 -\frac{1}{\sqrt{2\pi}}\int\limits_{-\infty}^v  e^{\frac{-t^2}{2}} \d t\Big\vert \leq 
 \frac{K}{\sqrt{n}\sigma a'}+\frac{1}{\pi} \frac{C'(1+ \sqrt{n}\sigma a'
 \rho_1^{n-1})}{\sqrt{n}}\leq  \frac{C''}{\sqrt{n}}.
 \end{equation*}
 This completes the proof of the theorem.
 \end{proof}


\section{Spectral decomposition of $\Lambda_h(it)$ for $t$ real}

We begin with a version of the classical 
Ionescu-Tulcea and Marinescu theorem  (see \cite{tm,no}).

\begin{thm}\label{thm_TM} Let $\mathcal{V}, $ $\mathcal{L}$ and $\mathcal{K}$ be three
Banach spaces respectively endowed with   
norms $ \Vert \cdot \Vert_{\mathcal{V}} ,$ $\Vert\cdot  \Vert_{\mathcal{L}}
,$   $\Vert\cdot  \Vert_{\mathcal{K}}$. Suppose that  
$\mathcal{V}\subset\mathcal{L}\subset\mathcal{K}$ 
and that 
there is a positive constant $C$ satisfying the following properties
\begin{itemize}
\item[(a)] If $(\phi_n)_{n=1}^{\infty}\subset \mathcal {V}$  
with  $ \Vert \phi_n\Vert_{\mathcal{V}}\leq 1$, then 
there is a subsequence $(\phi_{n_j})_{j=1}^{\infty} $ and  an element $\phi\in  \mathcal {V}$  such that
$ \lim\limits_{j\to\infty}\Vert \phi_{n_j}-\phi\Vert_{\mathcal{L}}=0 $ 
and $\Vert\phi\Vert_{\mathcal {V}}\leq C.$ 
\item[(a)']  $\Vert  \phi\Vert_{\mathcal{K}}
\leq C\Vert\phi\Vert_{\mathcal{L}},$  $\phi\in\mathcal{L}.$
\end{itemize}
Let $T :  \mathcal{K} \rightarrow \mathcal{K} $ be a bounded  linear operator such that 
$T(\mathcal{V})\subset\mathcal{V} .$ Assume that
\begin{itemize}
\item[(b)] (Doeblin--Fortet inequality) There are two 
positive constants $\alpha<1 $ and  $\beta$  such that
$  \Vert  T(\phi)\Vert_{\mathcal{V}}\leq  \alpha  \Vert   \phi
\Vert_{\mathcal{V}} +\beta \Vert   \phi \Vert_{\mathcal{L}},$  $\phi\in\mathcal{V}.$
\item[(c)]   $\rho(T|_{\mathcal{V}})\leq 1.$  
\item[(d)]   $T|_{\mathcal{V}}$ admits no eigenvalue of modulus $1,$  that is, if 
$\phi\in \mathcal{V}\setminus \{0\}$ and $\lambda\in\C$ satisfy
$T(\phi)=\lambda\phi,$ 
then $\vert\lambda\vert\not=1.$  
  \end{itemize}
Then $\rho(T|_{\mathcal{V}}) < 1.$
\end{thm}

\begin{proof}  It suffices  to prove for $\vert \lambda\vert=1$ that
 $(I-\lambda T)$ admits a bounded linear inverse from $\mathcal{V}$ onto
$\mathcal{V}.$  
Fix a $\lambda_0\in\C$ such that $\vert\lambda_0\vert  =1$ and an
arbitrary $\psi\in\mathcal{V}.$ 
We would like to find
a suitable $\phi\in\mathcal{V}$ such that $(I-\lambda_0 T)\phi=\psi.$ To this end
pick an arbitrary sequence  $ (\lambda_n)_{n=1}^{\infty}\subset \C$, 
$\vert\lambda_n\vert <1$, converging to $\lambda_0$.

Observe that by (c),
 for every $\lambda\in\C$ with $\vert\lambda\vert <1,$
$\rho(\lambda T)=\vert\lambda\vert \rho(T) <1.$  Therefore,
the series $(I-\lambda T)^{-1}:=\sum_{n\geq 0}\lambda^n T^n $ defines
 a bounded 
linear operator from $\mathcal{V}$ onto $\mathcal{V}$
which is  the inverse of $(I-\lambda T).$  Consequently, for every $n\geq 1,$ one may find
$\phi_n\in\mathcal{V}$ such that 
\begin{equation}\label{thm_TM_eq1}
(I-\lambda_n T)\phi_n=\psi.
\end{equation}
We claim that 
\begin{equation}\label{thm_TM_eq2}
 M:=\sup\limits_{n\geq 1}\Vert\phi_n\Vert_{\mathcal{L}}  <\infty.
\end{equation}
If not, there is a subsequence  $(\phi_{n_j})_{j=1}^{\infty}$ such that
$\lim_{j\to\infty} \Vert  \phi_{n_j}\Vert_{\mathcal{L}}=\infty.$ To simplify the notation assume that
$n_j=j$. Put  $\phi^{'}_j:=\frac{\phi_j}{\Vert \phi_j\Vert_{\mathcal{L}}}$. Then
(\ref{thm_TM_eq1}) may be rewritten as
\begin{equation} \label{thm_TM_eq3}
\phi_j^{'}=\lambda_j T\phi_j^{'}+\frac{\psi}{\Vert \phi_j\Vert_{\mathcal{L}}}.
\end{equation}
Consequently, using (b), we get
\begin{equation}  \label{thm_TM_eq4}
\Vert \phi_j^{'}\Vert_{\mathcal{V}}\leq \Vert T\phi_j^{'}\Vert_{\mathcal{V}}+\frac{\Vert \psi\Vert_{\mathcal{V}}}{\Vert \phi_j\Vert_{\mathcal{L}}}
\leq \alpha  \Vert \phi_j^{'}\Vert_{\mathcal{V}}+ \beta + \frac{\Vert \psi\Vert_{\mathcal{V}}}{\Vert \phi_j\Vert_{\mathcal{L}}}.
\end{equation}
Hence, 
\begin{equation}  \label{thm_TM_eq5}
\Vert \phi_j^{'}\Vert_{\mathcal{V}}\leq  \frac{ \beta + \frac{\Vert \psi\Vert_{\mathcal{V}}}{\Vert \phi_j\Vert_{\mathcal{L}}}}{1-\alpha}.
\end{equation}
Since by assumption $\lim_{j\to\infty} \Vert  \phi_{j}\Vert_{\mathcal{L}}=\infty,$ it follows from the latter estimate that 
 $(\phi_{j}^{'})_{j=1}^{\infty}$  is a bounded sequence in $\mathcal{V}.$  Therefore, by (a), there is a subsequence
 of $(\phi_{j}^{'})_{j=1}^{\infty}$  which converges  in $\mathcal{L}$ to an element $\phi^{'}\in\mathcal{V}.$
 Consequently,
 letting $j$ tend to $\infty$ in (\ref{thm_TM_eq3}) and using (a)' and 
the continuity of $T :\mathcal{K}\rightarrow\mathcal{K}$,  we obtain  
 \begin{equation*}
 \phi^{'}=\lambda_0 T\phi^{'}\qquad \text{and}\qquad \Vert \phi^{'}\Vert_{\mathcal{L}}=1,
 \end{equation*} 
 which contradicts (d). Hence, the proof of (\ref{thm_TM_eq2}) is complete.
 
 Using (\ref{thm_TM_eq2}) and arguing as in (\ref{thm_TM_eq4})--(\ref{thm_TM_eq5}), we can show that
 \begin{equation} \label{thm_TM_eq6}   
\Vert \phi_n\Vert_{\mathcal{V}}\leq  \frac{ M\beta +  \Vert \psi\Vert_{\mathcal{V}} }{1-\alpha}.
\end{equation}
By (a),  there is a subsequence
 of $(\phi_{n} )_{n=1}^{\infty}$  which converges  in $\mathcal{L}$ to
 an element $\phi \in\mathcal{V}.$ This, combined with
 (\ref{thm_TM_eq1}), implies that  $(I-\lambda_0 T)\phi=\psi.$ Hence,
 $(I-\lambda_0 T)$  
is onto  $\mathcal{V}$.  Recall from (c)--(d) that 
$(I-\lambda_0 T)$ is one-to-one. 
If $(\psi^{(n)})\subset \mathcal{V}$ and $(\phi^{(n)})\subset
 \mathcal{V}$ satisfy $\Vert \psi^{(n)}\Vert_{\mathcal{V}}= 1$ and $(I-\lambda_0
 T)\phi^{(n)}=\psi^{(n)}$, arguing as in    (\ref{thm_TM_eq4})--(\ref{thm_TM_eq6})
we show that $(\phi^{(n)})$ is bounded in $\mathcal{L}$, and then in
 $\mathcal{V}$. 
Therefore, $(I-\lambda_0
T)^{-1}$  exists and is bounded. This completes the proof.
\end{proof}

The remaining of this section is devoted to an application of Theorem
\ref{thm_TM} in the case  $ \mathcal{K}  :=L^1$,
$\mathcal{L}  :=L^{\frac{4+2\epsilon}{\epsilon}}$ and $\mathcal{V} :=W^{1,2}.$

\begin{prop}\label{prop_DF_W}
Let $h$ be  a real-valued function in $W^{1,2+\epsilon} $, $\epsilon>0.$ 
 Then, the hypotheses  (a)--(b) of
 Theorem \ref{thm_TM} are fulfilled with
 $ \mathcal{K}  :=L^1$,
$\mathcal{L}  :=L^{\frac{4+2\epsilon}{\epsilon}}$,
$\mathcal{V}:=W^{1,2}$ and 
 $T:=\Lambda_h(it )$,   $t \in\R$.
 \end{prop}
\begin{proof} The hypothesis (a) is a consequence of the Sobolev embedding
  theorem. The hypothesis (a)' is obvious. The operator
$\Lambda_h(it ):\K\rightarrow \K$ is  bounded since 
 for every $\phi\in L^1 ,$ we have 
\begin{equation}\label{eq1_prop_DF_W}
\begin{split}
\Vert  \Lambda_h(it )(\phi)\Vert_{L^1 } &=\left\Vert  
\Lambda\left(e^{i t  h}\phi\right)\right\Vert_{L^1 } 
\leq \Vert  \Lambda\left( \vert\phi\vert\right)\Vert_{L^1 }\\
&=\int\limits_{\P^1}\Lambda\left( \vert\phi\vert\right)\omega
={1\over d} \int\limits_{\P^1}\vert\phi\vert f^{\ast}(\omega)\leq C\Vert \phi\Vert_{L^1 }.
\end{split}
\end{equation}

We now check condition (b).
Let $\phi\in W^{1,2} .$ 
We first estimate $\Vert \d \Lambda_h(it )(\phi)\Vert_{L^2 }.$
To this end we  repeatedly apply the Cauchy-Schwarz inequality.
This inequality implies, in particular, that $if_*(\alpha)\wedge f_*(\overline \alpha)\leq
df_*(i\alpha\wedge\overline\alpha)$ or equivalently
$i\Lambda(\alpha)\wedge\Lambda(\overline \alpha) \leq
\Lambda(i\alpha\wedge\overline\alpha)$, where $\alpha$ is a $(1,0)$-form.
Recall that $\d=\partial+\overline\partial$.
Using Proposition \ref{prop_pertubated_oper}(4) we obtain
\begin{eqnarray*}
\lefteqn{i\partial \Lambda_h(it )(\phi)\wedge \overline{\partial
    \Lambda_h(it )(\phi)}} \\
&=&\frac{i}{d^{2}} f_{\ast}\big(  \partial (e^{it  h}\phi)\big) \wedge 
f_{\ast}\big(  \overline{\partial} (e^{-it  h}\overline{\phi})\big)\\
&\leq & \frac{1}{d}  f_{\ast}\big( i \partial (e^{it   h}\phi)  \wedge 
  \overline{\partial} (e^{-it  h}\overline{\phi})\big)\\
  &=& \frac{1}{d}  f_{\ast}\big( i \partial  \phi \wedge 
  \overline{\partial \phi}\big)+
  \frac{1}{d}  f_{\ast}\big(i \phi e^{-it  h} \partial e^{it   h}  \wedge 
  \overline{\partial\phi}  \big)\\
  && +\frac{1}{d}  f_{\ast}\big( i \overline{\phi}
  e^{it  h} \partial \phi\wedge \overline{\partial}e^{-it 
    h}    \big)   +\frac{1}{d}  f_{\ast}\big( i \vert\phi\vert^2
\partial e^{it  h}
\wedge \overline{\partial} e^{-it  h} \big)\\
  &\leq& \frac{11}{9d}  f_{\ast}\big( i \partial  \phi \wedge 
  \overline{\partial \phi}\big)+\frac{19}{d}
  f_{\ast}\big( i\vert\phi\vert^2 
\partial e^{it  h}\wedge \overline{\partial} e^{-it  h} \big).
\end{eqnarray*}
In the last line, we use the inequality $ab\leq {1\over 9}a^2+9b^2$.
 Consequently, since $\Vert \d\psi \Vert_{L^2}=2\Vert \partial \psi
 \Vert_{L^2}$ for every function $\psi$ and
 $\int_{\P^1}f_*(\beta)=\int_{\P^1} \beta$ for every 2-form $\beta$,
 we have 
 \begin{eqnarray*}
\Vert \d \Lambda_h(it )(\phi)\Vert_{L^2 }^2&=&
2\int\limits_{\P^1} i\partial \Lambda_h(it )(\phi)\wedge 
\overline{\partial \Lambda_h(it )(\phi)}\\
  &\leq& \frac{22}{9d} \int\limits_{\P^1} f_{\ast}\big( i \partial  \phi \wedge 
  \overline{\partial \phi}\big)+\frac{38}{d} 
\int\limits_{\P^1} (f^n)_{\ast}\big( i\vert\phi\vert^2 \partial 
e^{it  h}\wedge \overline{\partial} 
  e^{-it  h}  \big)\\
 & \leq &\frac{11}{18}\Vert \d  \phi\Vert_{L^2}^2+ 19
 \int\limits_{\P^1}  i \vert\phi\vert^2 \partial e^{it  h}\wedge \overline{\partial} 
  e^{-it  h}   .
\end{eqnarray*}
Applying  Cauchy--Schwarz inequality and   H{\"o}lder's inequality, we may estimate the  last integral    
$$\int\limits_{\P^1}  i \vert\phi\vert^2 \partial e^{it  h}\wedge \overline{\partial} 
 e^{-it  h}   =  \vert t \vert^2
 \int\limits_{\P^1} i  \vert\phi\vert^2 \partial   
h\wedge \overline{\partial}   h    
  \leq  \vert t \vert^2
\Vert \phi\Vert_{L^{\frac{4+2\epsilon}{\epsilon}}}^{ 2}
   \Vert \d h\Vert_{L^{2+\epsilon}}^{ 2}.$$
Finally, we have shown that
\begin{equation}\label{eq2_prop_DF_W}
\Vert \d \Lambda_h(it )(\phi)\Vert_{L^2}\leq 
\sqrt{\frac{11}{18}}\Vert \d  \phi \Vert_{L^2} +
C(h,\epsilon,t )  \Vert \phi\Vert_{ L^{\frac{4+2\epsilon}{\epsilon}}} .
\end{equation}
On the other hand, we  deduce from (\ref{PS2_inequality}) and
 (\ref{DS_inequality_W_1})--(\ref{DS_inequality_W_2}) 
   that
\begin{equation*}
\Vert   \Lambda_h(it )(\phi)\Vert_{L^2}\leq  \Vert   \Lambda (\vert\phi\vert)\Vert_{L^2}
\leq \int\Lambda\big(|\phi|^2\big)\omega = {1\over d}\int |\phi|^2 f^*(\omega) 
\leq  C\Vert\phi\Vert_{L^2}\leq  C'\Vert\phi\Vert_{L^{4+2\epsilon\over
  \epsilon}}.
\end{equation*}
  This, combined with   (\ref{eq2_prop_DF_W}), implies 
the Doeblin--Fortet inequality with $\alpha:=\sqrt{\frac{11}{18}},$
  and  $\beta:= C' + C(h,\epsilon,t )$.
\end{proof}

\begin{prop}\label{prop_spec_rad_W}
Under the hypothesis
 of Proposition \ref{prop_DF_W}, the hypothesis (c) in Theorem \ref{thm_TM}
 holds, i.e. we have  
$\rho\big(\Lambda_h(it )|_{W^{1,2}}\big)\leq 1$  for every
 $t\in\R$.
 \end{prop}
\begin{proof}
First, using Proposition \ref{prop_pertubated_oper}(4) and
(\ref{PS2_inequality})-(\ref{DS_inequality_W_2}),
we have
\begin{eqnarray} \label{eq_module_phi}
\Vert   \Lambda_h^n(it )(\phi)\Vert_{L^2}
& = & \Vert \Lambda^n(e^{itS_nh}\phi)\Vert_{L^2} \leq
\Vert\Lambda^n(\vert\phi\vert)\Vert_{L^2} \nonumber \\
& \leq & \Vert\Lambda^n(\vert\phi\vert)\Vert_{W^{1,2}} \leq
C\Vert   \vert\phi\vert  \Vert_{W^{1,2}} \leq C'\Vert   \phi  \Vert_{W^{1,2}}.
\end{eqnarray}

Now, we use Proposition \ref{prop_pertubated_oper}(4) and 
the Cauchy-Schwarz inequality repeatedly to estimate $\Vert \d
\Lambda_h^n(it )(\phi)\Vert_{L^2}.$  
As in Proposition \ref{prop_DF_W}, we obtain
\begin{eqnarray}\label{eq1_prop_spec_rad_W}
\lefteqn{i\partial \Lambda_h^n(it )(\phi)\wedge \overline{\partial
    \Lambda_h^n(it )(\phi)}}\nonumber\\
&= & i\Lambda^n\big(  \partial (e^{it  S_n h}\phi)\big)   
\wedge  \Lambda^n\big(  \overline{\partial} (e^{-it  S_n
  h}\overline{\phi})\big) \nonumber\\
 &\leq & i \Lambda^n\big(  it   e^{it  S_nh}\phi \partial  
 S_n h+ e^{it  S_n h} \partial\phi\big)
  \wedge \Lambda^n\big(-it   e^{-it  S_nh}\overline{\phi} 
\!\ \overline{\partial} S_n h  +  e^{-it  S_n h} \overline{\partial\phi}\big) \nonumber  \\
&\leq &  2 i \Lambda^n\big(   e^{it  S_nh}  \partial   \phi\big)
 \wedge \Lambda^n\big(    e^{-it  S_nh}  \overline{\partial
   \phi}\big) \nonumber\\
&& +  2t ^2 i\Lambda^n\big(    e^{it  S_nh}\phi \partial   S_n h \big)
 \wedge \Lambda^n\big(   e^{-it  S_nh}\overline{\phi}\!\ 
 \overline{\partial} S_n h  \big) \nonumber\\
 &  \equiv & I+ 2t ^2 II  .
\end{eqnarray}
For the first term we have
 \begin{equation}\label{eq2_prop_spec_rad_W}
 I\leq  2\Lambda^n \big(i\partial   \phi 
 \wedge  \overline{\partial \phi}\big).
 \end{equation}
For the second one we have
 \begin{eqnarray*}
 II&=& i \Lambda^n\Big(    e^{it  S_nh}\phi
 \sum\limits_{k=0}^{n-1}  \partial  (h\circ f^k) \Big)
 \wedge \Lambda^n\Big(   e^{-it 
 S_nh}\overline{\phi}\sum\limits_{k=0}^{n-1} 
\overline{\partial} (h\circ f^k)  \Big)\\
 &\leq& n  \sum\limits_{k=0}^{n-1} i \Lambda^n\big(    e^{it 
 S_nh}\phi  \partial   (h\circ f^k) \big)
 \wedge \Lambda^n\big(   e^{-it  S_nh}\overline{\phi} \!\ 
 \overline{\partial} (h\circ f^k)  \big)\\
  &=& n  \sum\limits_{k=0}^{n-1} i \Lambda^{n-k}\big(\Lambda^k\big(    
e^{it  S_nh}\phi  \partial   (h\circ f^k) \big)\big)
 \wedge \Lambda^{n-k}\big(\Lambda^k\big(   e^{-it 
 S_nh}\overline{\phi}\!\  \overline{\partial} (h\circ f^k) \big) \big)\\
 &\leq& n  \sum\limits_{k=0}^{n-1} \Lambda^{n-k}\Big(i\Lambda^k\big(    
e^{it  S_nh}\phi  \partial   (h\circ f^k) \big) 
 \wedge  \Lambda^k\big(   e^{-it  S_nh}\overline{\phi}\!\ 
 \overline{\partial} (h\circ f^k) \big) \Big).
\end{eqnarray*}
Observe that  for any complex-valued function $\psi$, we have 
\begin{equation*}
i\Lambda^k\big(\psi \partial   (h\circ f^k) \big) \wedge
\Lambda^k\big(\psi\overline{ \partial}   (h\circ f^k) \big)
= i\Lambda^k (\psi) \partial h\wedge \Lambda^k(\psi)\overline{\partial} h
\leq  i\Lambda^k (\vert\psi\vert)^2  \partial h
\wedge\overline{\partial} h.
\end{equation*}
Then
\begin{equation*}
 II \leq n  \sum\limits_{k=0}^{n-1} \Lambda^{n-k}\big( i\Lambda^k(   \vert\phi\vert)^2  \partial    h   
 \wedge    \overline{\partial} h  \big).
\end{equation*}
 This, combined with  (\ref{eq1_prop_spec_rad_W})--(\ref{eq2_prop_spec_rad_W}), implies that  
\begin{eqnarray*}
\Vert \d \Lambda_h^n(it )(\phi)\Vert_{L^2}^2 &\leq&
\frac{4}{d^n} \int\limits_{\P^1} (f^n)_{\ast} \big(i\partial   \phi 
 \wedge  \overline{\partial \phi}\big) \\
&& +  4nt ^2  \sum\limits_{k=0}^{n-1}
 \frac{1}{d^{n-k}}\int\limits_{\P^1} (f^{n-k})_{\ast} 
\big( i \Lambda^k(   \vert\phi\vert)^2  \partial    h   
 \wedge    \overline{\partial} h   \big)\\
 &=& \frac{4}{d^n} \int\limits_{\P^1}  i\partial   \phi 
 \wedge  \overline{\partial \phi} +
 4nt ^2  \sum\limits_{k=0}^{n-1} \frac{1}{d^{n-k}} 
\int\limits_{\P^1}  i\Lambda^k(   \vert\phi\vert)^2  \partial    h   
 \wedge    \overline{\partial} h .   
 \end{eqnarray*}
 Applying H{\"o}lder's inequality,  (\ref{PS2_inequality}) and 
(\ref{eq_module_phi})  yields that 
 \begin{equation*}
\int\limits_{\P^1}  i\Lambda^k(   \vert\phi\vert)^2  \partial    h   
 \wedge    \overline{\partial} h \leq    \Vert
 \Lambda^k(\vert\phi\vert) 
\Vert_{L^{\frac{4+2\epsilon}{\epsilon}}}^{ 2}
   \Vert \d h \Vert_{L^{2+\epsilon}}^{ 2}
   \leq C_{\epsilon} \Vert  \phi  \Vert_{W^{1,2}}^{2}
   \Vert \d h \Vert_{L^{2+\epsilon}}^{ 2}.
 \end{equation*}
 Since $\sum \frac{1}{d^{n-k}}\leq 2,$ it follows that
 \begin{equation}\label{eq_last_prop_spec_rad_W} 
\Vert \d \Lambda_h^n(it )(\phi)\Vert_{L^2}\leq\Big(
\sqrt{\frac{4}{d^n}}  +  \sqrt{ 8nt ^2   C_{\epsilon}
  } 
  \Vert \d h \Vert_{L^{2+\epsilon}} \Big)\Vert  \phi  \Vert_{W^{1,2}}   .   
 \end{equation}
Combining  this and  (\ref{eq_module_phi}), the desired conclusion follows. 
\end{proof}


\section{Local Central Limit Theorem}

First we introduce the following notion. 
\begin{defi}\label{Condition_H} 
A real-valued function $h$ is a multiplicative  $W^{1,2}$-cocycle if
there are $t>0$, $s\in\R$ and  $\phi\in W^{1,2}$, $\phi$ non zero in $L^1(\mu)$, such that 
$ e^{it h(z)}\phi(z)=e^{is}\phi(f(z))$ $\mu$-a.e.
\end{defi}

We have the following proposition.
\begin{prop}\label{sufficiency}
Let $h$ be a real-valued function in $ W^{1,2+\epsilon}.$ 
Then $h$ is not a multiplicative  $W^{1,2}$-cocycle if and only if 
the spectral radius of  $\Lambda_h(it)  :W^{1,2}\rightarrow W^{1,2}$
is strictly smaller than $1$ for every $t>0$.
\end{prop} 
\begin{proof}
First assume that $h$ is not a  multiplicative  $W^{1,2}$-cocycle.
Suppose in order to get a contradiction that
there is $t>0$ such that
$\rho\big(\Lambda_h(it)\vert_{W^{1,2}}\big)\geq 1$.
By Propositions \ref{prop_DF_W}, \ref{prop_spec_rad_W} and  Theorem \ref{thm_TM}, 
there are $\lambda\in \C$ and   $\phi\in W^{1,2}\setminus\{0\}$  such that
$\Lambda_h(it) (\phi)=\lambda \phi$  and $\vert\lambda\vert=1.$
It follows that $\Lambda(e^{ith}\phi)=\lambda\phi$, then $\Lambda(\vert\phi\vert)\geq\vert\phi\vert$. Since
$\mu(\Lambda(\vert\phi\vert))=\mu(\vert\phi\vert)$, 
we deduce that $\Lambda(\vert\phi\vert)= \vert\phi\vert$ $\mu$-a.e.
We also obtain $\vert\phi\vert=\Lambda^n(\vert\phi\vert)$, and since
$\Lambda^n(\vert\phi\vert)$ converges in $L^1(\mu)$ to a constant,
$|\phi|$ is constant $\mu$-a.e.
 Next, we rewrite   $\Lambda(e^{ith}\phi)=\lambda \phi$ as
\begin{equation*}
          \frac{1}{d}\sum\limits_{w\in f^{-1}(z)} (e^{ith}\phi)(w)
=\lambda\phi(z)        ,\qquad z\in\P^1.
\end{equation*}
Since $|\phi|$ is constant, this is only possible if 
\begin{equation*}
e^{it h(z)}\phi(z)= \lambda\phi(f(z)),  \qquad \mu  \text{-a.e.}\ z\in \P^1.
\end{equation*}
Since $h$ is not a cocycle, this implies that     $\phi=0$ for $\mu$-a.e. $z\in \P^1.$
  Next, consider the function  $\Lambda^n(\vert\phi\vert),$ $n\in\N.$
  Since one has just shown that  $\mu(\vert\phi\vert)=0,$  it follows
  from (\ref{DS_inequality_W_2}) and the identity $\Lambda_h(it)(\phi)=\lambda\phi$
that
  \begin{equation*}
  \Vert\phi\Vert_{L^1}=\left\Vert  \Lambda_h^n(it)(\phi)\right\Vert_{L^1}
  =\left\Vert  \Lambda^n\big( e^{itS_nh}\phi\big )\right\Vert_{L^1}
\leq  \left\Vert  \Lambda^n ( \vert \phi\vert  )\right\Vert_{L^1}\leq 
\frac{C}{d^{n/2}} \Vert \phi\Vert_{W^{1,2}}. 
  \end{equation*} 
Letting $n$ tend to $\infty,$ we  obtain $\phi=0,$ which is a contradiction.

Now, assume that $\rho\big(\Lambda_h(it)\vert_{W^{1,2}}\big)<1$. Then
for $\phi\in W^{1,2}$ we have $\Lambda_h^n(it)(\phi)\rightarrow 0$ in
$W^{1,2}$. It follows that $\Lambda_h^n(it)(\phi)\rightarrow 0$
$\mu$-a.e. If $e^{ith(z)}\phi(z)=e^{is}\phi(f(z))$ $\mu$-a.e. then
$\Lambda_h^n(it)(\phi)=e^{ins}\phi$ $\mu$-a.e. Hence $\phi=0$
$\mu$-a.e. This implies that $h$ is not a  multiplicative  $W^{1,2}$-cocycle.
\end{proof}

Now we are able to state the last main result.
\renewcommand{\thethmspec}{Theorem C}
  \begin{thmspec} 
  We keep the hypothesis and the notation of Theorem A.  Assume that $h$
 is not a  multiplicative  $W^{1,2}$-cocycle.
   Then for every bounded interval $\Delta\subset \R,$ the following convergence
  holds uniformly in $x\in\R$
  \begin{equation*}
 \lim\limits_{n\to\infty}\Big\vert\sigma\sqrt{n}\nu\big\lbrace
   x+S_n h  \in\Delta \big\rbrace
 -\frac{1}{\sqrt{2\pi}}   e^{-\frac{x^2}{2\sigma^2n} }m(\Delta)\Big\vert=0 ,
 \end{equation*}
 where $m(\Delta)$ denotes the length of $\Delta.$
 \end{thmspec}
 \begin{proof}
 We follow the proof given by Breiman (see \cite[pp. 224-227]{bl}) in the context of independant
 random variables.  The proof of the above equality is  reduced   to showing that
 \begin{equation}\label{eq1_TheoremC}  
 \lim\limits_{n\to\infty}\Big\vert\sigma\sqrt{n}\int\limits_{\P^1}\phi(x+S_nh)\d\nu 
 -\frac{1}{\sqrt{2\pi}} \int\limits_{-\infty}^{\infty}\phi(t)
 e^{-\frac{x^2}{2\sigma^2n} }\d t
\Big\vert \equiv\lim\limits_{n\to\infty}
 \Big\vert  \frac{A_n(x)}{\sqrt{2\pi}}\Big\vert=0
 \end{equation}
 for every  real-valued function $\phi\in L^1(\R)$ whose Fourier transform
 $\widehat{\phi}(x):=\frac{1}{\sqrt{2\pi}}\int_\R \phi(t)e^{-itx} \d t $   
is  a  continuous function supported 
in an  interval $[-\delta,\delta]$, $\delta>0.$    
 
 Observe that 
\begin{eqnarray*}  
  \sigma\sqrt{n}\int\limits_{\P^1}\phi(x+S_nh)\d\nu 
  &=&\frac{\sigma\sqrt{n}}{\sqrt{2\pi}} \int\limits_{\P^1}
\Big(  \int\limits_{-\delta}^{\delta}\widehat{\phi}(t)e^{it(x+S_nh)} \d t\Big) g\d \mu\\
 &=& \frac{\sigma\sqrt{n}}{\sqrt{2\pi}} \int\limits_{-\delta}^{\delta}
  \widehat{\phi}(t)e^{itx }
    \Big(  \int\limits_{\P^1} e^{it S_nh}  g\d\mu \Big) \d t\\
  &=& \frac{ 1}{\sqrt{2\pi}}
  \int\limits_{-\delta\sigma\sqrt{n}}^{\delta\sigma\sqrt{n}} 
\widehat{\phi}\Big(\frac{t}{\sigma\sqrt{n}}\Big)
  e^{\frac{itx}{\sigma\sqrt{n}} }    
\Big(  \int\limits_{\P^1} e^{it \frac{S_nh}{\sigma\sqrt{n}}}  g\d\mu \Big) \d t\\
 &=& \frac{ 1}{\sqrt{2\pi}}
  \int\limits_{-\delta\sigma\sqrt{n}}^{\delta\sigma\sqrt{n}} 
\widehat{\phi}\Big(\frac{t}{\sigma\sqrt{n}}\Big)
  e^{\frac{itx}{\sigma\sqrt{n}} }    
\Big(  \int\limits_{\P^1}  \Lambda_h^n\big(\frac{it}{\sigma\sqrt{n}}\big)(  g)\d\mu \Big) \d t 
  .\end{eqnarray*}
We also need the following identities, 
\begin{equation*}
\widehat{\phi}(0)=\frac{1}{\sqrt{2\pi}}
 \int\limits_{-\infty}^{\infty}\phi(t) \d t \quad\text{and}\quad
 e^{-\frac{x^2}{2\sigma^2n} }=
\frac{1}{\sqrt{2\pi}} \int\limits_{-\infty}^{\infty}   e^{it\frac{x}{\sigma\sqrt{n}} }e^{-\frac{t^2}{2}}\d t.
\end{equation*}
If we replace in (\ref{eq1_TheoremC}) we obtain
\begin{equation*}
 A_n(x)= \int\limits_{-\delta\sigma\sqrt{n}}^{\delta\sigma\sqrt{n}}
 \widehat{\phi}\Big(\frac{t}{\sigma\sqrt{n}}\Big)
  e^{\frac{itx}{\sigma\sqrt{n}} }    \Big(  
\int\limits_{\P^1}  \Lambda_h^n\big(\frac{it}{\sigma\sqrt{n}}\big)(  g)\d\mu \Big) \d t -
 \widehat{\phi}(0)  \int\limits_{-\infty}^{\infty}   
e^{it\frac{x}{\sigma\sqrt{n}} }e^{-\frac{t^2}{2}}\d t.
\end{equation*}
In order to estimate $A_n(x),$ we divide it into three pieces 
\begin{equation*}
A_n(x)\equiv A_n^1(x)+A_n^2(x)+A_n^3(x).
\end{equation*}
 They are given   by the following formulas where the constant
 $\alpha\in (0,\delta)$  will be  determined later on : 
\begin{equation}\label{eq2_TheoremC}
\begin{split}
A_n^1(x)&:= \int\limits_{-\alpha\sigma\sqrt{n}}^{\alpha\sigma\sqrt{n}} 
\widehat{\phi}\Big(\frac{t}{\sigma\sqrt{n}}\Big)
  e^{\frac{itx}{\sigma\sqrt{n}} }    
\Big(  \int\limits_{\P^1}  \Lambda_h^n\big(\frac{it}{\sigma\sqrt{n}}\big)(  g)\d\mu \Big) \d t -
 \widehat{\phi}(0)  \int\limits_{-
   \alpha\sigma\sqrt{n}}^{\alpha\sigma\sqrt{n}}   
e^{it\frac{x}{\sigma\sqrt{n}} }e^{-\frac{t^2}{2}}\d t,\\
 A_n^2(x)&:=\int\limits_{\alpha<\vert\frac{t}{\sigma\sqrt{n}}
   \vert<\delta}  
\widehat{\phi}\Big(\frac{t}{\sigma\sqrt{n}}\Big)
  e^{\frac{itx}{\sigma\sqrt{n}} }    \Big(  
\int\limits_{\P^1}  \Lambda_h^n\big(\frac{it}{\sigma\sqrt{n}}\big)(  g)\d\mu \Big) \d t,\\
 A_n^3(x)&:= \widehat{\phi}(0)  
\int\limits_{\alpha<\vert\frac{t}{\sigma\sqrt{n}}\vert }   
e^{it\frac{x}{\sigma\sqrt{n}} }e^{-\frac{t^2}{2}}\d t.
 \end{split}
\end{equation}
To estimate  $A_n^1(x)$  we recall from  Theorem   \ref{thm_pertubation}  that
    for $\big\vert \frac{t}{\sigma\sqrt{n}}\big\vert <a $ 
\begin{equation*}
\int\limits_{\P^1}\Lambda_h^n\Big(
\frac{t}{\sigma\sqrt{n}}\Big)(g)\d\mu
=\lambda^n\Big( \frac{t}{\sigma\sqrt{n}}\Big) 
\int\limits_{\P^1}\Phi\Big( \frac{t}{\sigma\sqrt{n}}\Big)(g)\d\mu  +
\int\limits_{\P^1}\Psi^n \Big( \frac{t}{\sigma\sqrt{n}}\Big)(g)\d\mu .
\end{equation*}
On the other hand, we know from the proof of Theorem A that
\begin{equation*}
 \lim\limits_{n\to\infty} \lambda^n\Big( \frac{t}{\sigma\sqrt{n}}\Big) =e^{-\frac{t^2}{2}},\quad
\lim\limits_{n\to\infty}\int\limits_{\P^1}\Phi\Big( \frac{t}{\sigma\sqrt{n}}\Big)(g)\d\mu =1,\quad 
 \lim\limits_{n\to\infty} \int\limits_{\P^1}\Psi^n \Big( \frac{t}{\sigma\sqrt{n}}\Big)(g)\d\mu  =0
 \end{equation*}
and there is  a constant  $\alpha>0$  (we can choose $\alpha<\delta$) such that
for $\big\vert \frac{t}{\sigma\sqrt{n}}\big\vert <\alpha, $ 
\begin{equation*}\label{eq4_TheoremC}
\begin{split}
\Big\vert \lambda^n\Big( \frac{t}{\sigma\sqrt{n}}\Big) 
\int\limits_{\P^1}\Phi\Big(
\frac{t}{\sigma\sqrt{n}}\Big)(g)\d\mu\Big\vert &<  
C e^{-\frac{t^2}{4}},\\
 \Big\vert\int\limits_{\P^1}\Psi^n \Big(
   \frac{t}{\sigma\sqrt{n}}\Big)(g)\d\mu \Big\vert 
&< C\rho_1^n<C e^{-\frac{t^2}{4}}.
\end{split}
\end{equation*}
Then, 
\begin{equation*}
 \Big\vert  \widehat{\phi}\Big(\frac{t}{\sigma\sqrt{n}}\Big)
       \int\limits_{\P^1}  \Lambda_h^n\Big(\frac{it}{\sigma\sqrt{n}}\Big)(  g)\d\mu   -
 \widehat{\phi}(0)    e^{-\frac{t^2}{2}} \Big\vert \leq  C  
e^{-\frac{t^2}{4}}\Vert \widehat{\phi}\Vert_{L^{\infty}(\R)}, \qquad
       \frac{t}{\sigma\sqrt{n}}\in [-\alpha,\alpha]. 
\end{equation*}
  Observe that the right hand side of the latter estimate is
an integrable function.    
On the other hand, since $\widehat{\phi}$ is continuous,  we
deduce from the previous identities that 
\begin{equation*}
 \lim\limits_{n\to\infty} \widehat{\phi}\Big(\frac{t}{\sigma\sqrt{n}}\Big)
         \int\limits_{\P^1}  \Lambda_h^n\Big(\frac{it}{\sigma\sqrt{n}}\Big)(  g)\d\mu   -
 \widehat{\phi}(0)    e^{-\frac{t^2}{2}} =0, \qquad t\in\R.
\end{equation*}
Then the Lebesgue's  
dominated convergence theorem implies
\begin{equation*}
\lim\limits_{n\to\infty} \vert A_n^1(x)\vert=0
\end{equation*} 
uniformly in $x\in\R.$

Now  we  estimate   $A_n^2(x).$  For every $\tau\in [\alpha,\delta],$ %
  consider the operator   $T:=\Lambda_{h}(i\tau)\ : W^{1,2}\rightarrow W^{1,2}$.
Recall from Proposition 
\ref{sufficiency} that
$\rho\left(\Lambda_h(i\tau)\right) <1$.
Using Proposition \ref{prop_pertubated_oper} and a compactness argument 
 we may find  $C>0$ and  $0<\rho<1$ such that
 $\Vert\Lambda_h^n(i\tau)\Vert\leq C\rho^n$, $\tau\in[\alpha,\delta]$. Hence
 \begin{equation*}
 \Big\vert \int\limits_{\P^1}
 \Lambda_h^n\Big(\frac{it}{\sigma\sqrt{n}}\Big)(  g)\d\mu 
\Big\vert  \leq C' \rho^n
 \end{equation*}
 for $ \alpha\leq \big\vert\frac{t}{\sigma\sqrt{n}} \big\vert \leq \delta$.
 Consequently, we deduce from (\ref{eq2_TheoremC}) that
  \begin{equation*}
\lim\limits_{n\to\infty} \vert A_n^2(x)\vert\leq \lim\limits_{n\to\infty} C'' \sqrt{n}\rho^n=0
\end{equation*} 
uniformly in $x\in\R.$

 To estimate    $A_n^3(x)$  it suffices to see  that
 \begin{equation*}
\lim\limits_{n\to\infty} \vert A_n^3(x)\vert\leq \lim\limits_{n\to\infty}
 \int\limits_{\sqrt{n}\sigma\alpha< \vert t  \vert   }  e^{-\frac{t^2}{2}}\d t=0.
 \end{equation*}
The estimates for $A_n^1(x),$ $A_n^2(x)$ and $A_n^3(x)$ imply
  (\ref{eq1_TheoremC}). 
\end{proof}

  \smallskip

\noindent{\bf Acknowledgment.}     
The second named   author  wishes to express his gratitude to the Max-Planck Institut f{\"u}r Mathematik
in Bonn  for its hospitality and  support.


\end{document}